\documentclass[]{article}
\usepackage{amsfonts, amsmath, amssymb}
\usepackage{graphicx,color}
\usepackage{graphics}
\usepackage{psfig}

\def\N{\textrm{I\kern-0.21emN}}

\def\R{\textrm{I\kern-0.21emR}}
\def\Q{\textrm{l\kern-0.5emQ}}

\begin{document}

\title{Spurious caustics of \textit{Dispersion Relation Preserving} schemes}

\author{Claire David \footnotemark[1] , Pierre Sagaut \footnotemark[1]  \\
 \\
\footnotemark[1] Universit\'e Pierre et Marie Curie-Paris 6  \\
Institut Jean Le Rond d'Alembert, UMR CNRS 7190 \\
Bo\^ite courrier $n^0162$, 4 place Jussieu, F-75252 Paris
cedex 05, France\\
 }

\maketitle

\begin{abstract}

A linear dispersive mechanism leading to a burst in the $L_\infty$
norm of the error  in numerical simulation of polychromatic
solutions is identified. This local error pile-up corresponds to
the existence of spurious caustics, which are allowed by the
dispersive nature of the numerical error. From the mathematical
point of view, spurious caustics are related to extrema of the
numerical group velocity and are physically associated to
interactions between rays defined by the characteristic lines of
the discrete system. This paper extends our previous work about
classical schemes to dispersion-relation preserving schemes.

\end{abstract}


\section{Introduction}
\label{sec:intro}

\noindent The analysis and the control of numerical error in
discretized propagation-type equations is of major importance for
both theoretical analysis and practical applications. A huge
amount of works has been devoted to the analysis of the numerical
errors, its dynamics and its influence on the computed solution
(the reader is referred to classical books, among which
\cite{Hirsch,vich,lomax,tapan-book}). The emergence of
\textit{Dispersion-Relation-Preserving (DRP)} schemes \cite{Tam}),
which have the same dispersion relation as the original partial
difference equations, enables one to
 have very accurate high order finite difference schemes.\\

\noindent The two sources of numerical error are the dispersive
and dissipative properties of the numerical scheme, which are very
often investigated in unbounded or periodic domains thanks to a
spectral analysis.\\

\noindent It appears that existing works are mostly devoted to
linear, one-dimensional numerical models, such as the linear
advection equation

\begin{equation}
\label{transp}
 \frac{\partial u}{\partial t}+c \, \frac{\partial u}{\partial x}=0
\end{equation}

\noindent where $c$ is a constant uniform advection velocity.

\noindent The two sources of numerical error are the dispersive
and dissipative properties of the numerical scheme, which are very
often investigated in unbounded or periodic domains thanks to a
spectral analysis. Following this approach, a monochromatic wave
is used to measure the accuracy of the scheme. Such a tool is very
powerful and provides the user with a deep insight into the
discretization errors. But some results coming from practical
numerical experiments still remain unexplained, despite the linear
character of the discrete numerical model. As an example, let us
note the sudden growth of the numerical error for long range
propagation reported by Zingg \cite{zingg} for a large set of
numerical schemes, including optimized numerical schemes.

\noindent The usual modal analysis is almost always applied to
monochromatic reference solutions, with the purpose of analyzing
the error committed on both their amplitude and their phase,
leading to classical plots of the relative error as the function
of the Courant number and/or the number of grid points per
wavelength. Therefore, dispersive phenomena associated to
polychromatic solutions are usually not taken into account.\\

\noindent The present paper deals with the analysis of linear
dispersive mechanism which results in local error focusing, i.e.
to a sudden local error burst in the $L_\infty$ norm for
polychromatic solutions. This phenomena is reminiscent of the
physical one referred to as the caustic phenomenon in linear
dispersive physical models \cite{witham}, and will be referred to
as the spurious caustic phenomenon hereafter. It extends our
previous work \cite{Cl. David} to \textit{DRP} schemes. The
present analysis is restricted to interior stencil, and the
influence of boundary conditions will not be considered.

\noindent The paper is organized as follows. Main elements of caustic
theory of interest for the present analysis are recalled in
section \ref{Caustiques}. \textit{DRP} schemes are presented in section \ref{DRPStudy}. Their caustical analysis is exposed in section \ref{CaustiquesDRP}.  A numerical example is presented in section \ref{Example}.

\section{Caustics}
\label{Caustiques}

\noindent The solution of Eq. (\ref{transp}) is taken under the
form:

\begin{equation}
\label{depl} u(x,t,k)=e^{i\,(k\,x-\omega \, t)}
\end{equation}

\noindent where $\omega = {\xi}_\omega+i \,{\eta}_\omega $ is the
complex phase, and $k$ the real wave number. For dispersive waves,
it is recalled that  the group velocity $V_g (k)$ is defined as

\begin{equation}
\label{Vg} V_g (k) \equiv \frac{{\partial } \xi _\omega }{\partial
{k}}
\end{equation}

\noindent A caustic is defined as a focusing of different rays in
a single location. The equivalent condition is that the group
velocity exhibits an extremum, i.e. there exists at least one wave
number $k_c$ such that

\begin{equation}
\label{eq:caustic} \frac{{\partial } V_g}{\partial {k}}  (k_c) =0
\end{equation}

\noindent The corresponding physical interpretation is that wave
packets with characteristic wave numbers close to $k_c$ will
pile-up after a finite time and will remain superimposed for a
long time, resulting in the existence a region of high energy
followed by a region with very low fluctuation level.

\noindent The linear continuous model Eq. (\ref{transp}) is not
dispersive if the convection velocity $c$ is uniform, and
therefore the exact solution does not exhibits caustics since the
group velocity does not depend on $k$. The discrete solution
associated with a given numerical scheme will admit spurious
caustics, and therefore spurious local energy pile-up and local
sudden growth of the error, if the discrete dispersion relation is
such that the condition (\ref{eq:caustic}) is satisfied. For a
uniform scale-dependent convection velocity, such spurious
caustics can exist in polychromatic solutions only, since they are
associated to the superposition of wave packets with different
characteristic wave numbers.\\
\noindent Set:

\begin{equation}
k=\frac{\varphi \,\sigma}{c\,dt}
\end{equation}

\noindent The general dispersion relation associated with the
discrete scheme enables us to obtain the corresponding group velocity, given by:

\begin{equation}
\label{eq:Vg} V_g =h  \frac{{\partial} \xi _\omega }{\partial
{\varphi}}
  \end{equation}

\noindent The numerical solution will therefore admits spurious
caustics if

\begin{equation}
\frac{{\partial } V_g }{\partial {k}} = \frac{{\partial } V_g
}{\partial {\varphi}} \, \frac{{\partial } \varphi }{\partial {k}}=0
 \Longleftrightarrow
\frac{{\partial } V_g }{\partial {\varphi}}=0
          \end{equation}

\noindent The corresponding values of $\varphi$ and $k$ will be
respectively denoted ${\varphi}_c$ and $k_c$.

\noindent Spurious caustics are associated with characteristic
lines given by

\begin{equation}
 \frac{x}{t}=U_c
\end{equation}

\noindent  where

\begin{equation}
 U_c=V_g(\varphi_c)
\end{equation}

\section{\textit{DRP} schemes}
\label{DRPStudy}

\noindent The Burgers equation:

\begin{equation}
\label{Burgers}  u_t + c\, u\, u_x - \mu \,u_{xx} = 0,
\end{equation}

\noindent $c$, $\mu$ being real constants, plays a crucial role in
the history of wave equations. It was named after its use by
Burgers \cite{burger1} for studying
turbulence in 1939.\\

\noindent $i$, $n$ denoting natural integers, a linear finite
difference scheme for this equation can be written under the form:
\begin{equation} \label{scheme} \displaystyle \sum \alpha_{lm}\,u_{l}^{m}=0
              \end{equation}

\noindent where:
\begin{equation}
{u_l}^m=u\,(l\,h, m\,\tau)
\end{equation}
\noindent  $l\, \in \, \{i-1,\, i, \, i+1\}$, $m \, \in \,
\{n-1,\, n, \, n+1\}$, $j=0, \, ..., \, n_x$, $n=0, \, ..., \,
n_t$. The $ \alpha_{lm}$ are real coefficients, which depend on the mesh size $h$, and the time step $\tau$.\\
The Courant-Friedrichs-Lewy number ($cfl$) is defined as $\sigma = c \,\tau / h$ .\\
\noindent A numerical scheme is  specified by selecting
appropriate values of the coefficients $ \alpha_{lm}$. Then,
depending on them, one can obtain optimum schemes, for which the
error will be minimal.\\

 \noindent $m$ being a strictly positive integer, the first derivative $\frac{\partial
u}{\partial x}$ is approximated at the $l^{th}$ node of the
spatial mesh by:

\begin{equation}\label{approx}
 (\, \frac{\partial u}{\partial
x}\,)_l  \simeq
   \displaystyle \sum_{k=-m}^m \gamma_{k}\,u_{i+k}^n
\end{equation}
\noindent Following the method exposed by C. Tam and J. Webb in
\cite{Tam}, the coefficients $\\gamma_{k}$ are determined
requiring the Fourier Transform of the finite difference scheme
(\ref{approx}) to be a close approximation of the partial
derivative $ (\, \frac{\partial u}{\partial x}\,
)_l$.\\
\noindent (\ref{approx}) is a special case of:

\begin{equation}\label{approx_Cont}
 (\, \frac{\partial u}{\partial
x}\,)_l  \simeq \displaystyle \sum_{k=-m}^m \gamma_{k}\,u(x+k\,h)
\end{equation}

\noindent where $x$ is a continuous variable, and can be recovered
setting $x=l\,h$.\\
\noindent Denote by $\omega$ the phase. Applying the Fourier
transform, referred to by $\,\widehat{\, }$ , to both sides of
(\ref{approx_Cont}), yields:

\begin{equation}
\label{Wavenb}
 j\, \omega \, \widehat{u}  \simeq \displaystyle \sum_{k=-m}^m \gamma_{k}\,e^{\,j\,k\,\omega\,h}\, \widehat{u}
\end{equation}
\noindent  $j$ denoting the complex square root of $-1$.\\




\noindent Comparing the two sides of (\ref{Wavenb}) enables us to
identify the wavenumber $ \overline{\lambda}$ of the finite
difference scheme (\ref{approx}) and the quantity
$\frac{1}{j}\,{\displaystyle \sum_{k=-m}^m
\gamma_{k}\,e^{\,j\,k\,\omega\,h}}$, the
wavenumber of the finite difference scheme (\ref{approx}) is thus:

\begin{equation}
 \overline{\lambda}=-\,j\, \displaystyle \sum_{k=-m}^m \gamma_{k}\,e^{\,j\,k\,\omega\,h}
\end{equation}

\noindent To ensure that the Fourier transform of the finite
difference scheme is a good approximation of the partial
derivative $ (\, \frac{\partial u}{\partial x}\, )_l$ over the
range of waves with wavelength longer than $4\,h$, the a priori
unknowns coefficients $\gamma_{k}$ must be choosen so as to
minimize the integrated error:

\footnotesize
\begin{equation}\begin{array}{rcl}
 {\mathcal E} &=&\int_{-\frac{\pi}{2}}^{\frac{\pi}{2}} | \lambda \,h- \overline{\lambda}
 \,h|^2\,d(\lambda \,h)\\
 &=&\int_{-\frac{\pi}{2}}^{\frac{\pi}{2}} | \lambda \,h+j\, \displaystyle \sum_{k=-m}^m \gamma_{k}\,e^{\,j\,k\,\omega\,h} \,h|^2\,d(\lambda \,h)\\
 &=& \int_{-\frac{\pi}{2}}^{\frac{\pi}{2}} | \zeta+j\, \displaystyle \sum_{k=-m}^m \gamma_{k}\,\left  \lbrace \cos (\,k\,\zeta)+j\,\sin(\,k\,\zeta) \right \rbrace \,
 |^2\,d \zeta \\
 &=& \int_{-\frac{\pi}{2}}^{\frac{\pi}{2}} \left \lbrace \left [ \zeta- \displaystyle \sum_{k=-m}^m \gamma_{k} \,\sin(\,k\,\zeta) \right ]^2+
  \left [   \displaystyle \sum_{k=-m}^m \gamma_{k}\, \cos (\,k\,\zeta)  \right ]^2\, \right \rbrace  \,d
  \zeta \\
   &=& 2 \, \int_0^{\frac{\pi}{2}} \left \lbrace \left [ \zeta- \displaystyle \sum_{k=-m}^m \gamma_{k} \,\sin(\,k\,\zeta) \right ]^2+
  \left [   \displaystyle \sum_{k=-m}^m \gamma_{k}\, \cos (\,k\,\zeta)  \right ]^2\, \right \rbrace  \,d
  \zeta \\
      \end{array}
\end{equation}

\normalsize

\noindent The conditions that ${\mathcal E}$ is a minimum are:

\begin{equation}
 \frac {\partial{\mathcal E}}{\partial \gamma_i} =0 \,\,\,  , \,\,\,
i=-m, \ldots , \, m
\end{equation}

\noindent i. e.:

\begin{equation}
\label{RelDer}
  \int_0^{\frac{\pi}{2}} \left \lbrace
 \,-\,\zeta\,\sin(\,i\,\zeta)\, + \displaystyle \sum_{k=-m}^m \gamma_{k}
\,\cos\left(\,(k-i)\,\zeta \right)    \right \rbrace  \,d
  \zeta =0
\end{equation}

\noindent Changing $i$ into $-i$, and $k$ into $-k$ in the
summation yields:

\begin{equation}
  \int_0^{\frac{\pi}{2}} \left \lbrace
 \,\,\zeta\,\sin(\,i\,\zeta)\, + \displaystyle \sum_{k=-m}^m \gamma_{-k}
\,\cos\left(\,(-k+i)\,\zeta \right)    \right \rbrace  \,d
  \zeta =0
\end{equation}

\noindent i. e.:

\begin{equation}
  \int_0^{\frac{\pi}{2}} \left \lbrace
 \,\,\zeta\,\sin(\,i\,\zeta)\, + \displaystyle \sum_{k=-m}^m \gamma_{-k}
\,\cos\left(\,(k-i)\,\zeta \right)    \right \rbrace  \,d
  \zeta =0
\end{equation}

\noindent Thus:

\begin{equation}
\label{Int1}
  \int_0^{\frac{\pi}{2}}  \displaystyle \sum_{k=-m}^m \left \lbrace
  \gamma_{-k} + \gamma_{k}  \right \rbrace
\,\cos\left(\,(k-i)\,\zeta \right)    \,d
  \zeta =0
\end{equation}

\noindent which yields:

\begin{equation}
 \frac{\pi}{2} \,\left \lbrace
  \gamma_{-i} + \gamma_{i}  \right \rbrace +\displaystyle \sum_{k \neq  i,\, k=-m}^m
\left \lbrace \frac {
  \gamma_{-k} + \gamma_{k} }{k-i} \right \rbrace
\,\sin\left(\,(k-i)\,\frac{\pi}{2} \right)    \  =0
\end{equation}

\noindent which can be considered as a linear system of $2\,m+1$
equations, the unknowns of which are the $ \gamma_{-i} +
\gamma_{i} $, $i=-m,\,\ldots, \,m$. The determinant of this system
is not equal to zero, while it is the case of its second member:
the Cramer formulae give then, for $i=-m,\,\ldots, \,m$:

\begin{equation}
  \gamma_{-i} + \gamma_{i}    =0
\end{equation}

 \noindent or:

\begin{equation}
\label{RelCoeffDRP1}
  \gamma_{-i} =- \gamma_{i}
\end{equation}

 \noindent For $i = 0$, one of course obtains:

\begin{equation}
   \gamma_0=0
\end{equation}

 \noindent All this ensures:

\begin{equation}
\label{RelCoeffDRP} \displaystyle \sum _{k=-m}^m \gamma_{k}  =0
\end{equation}

 \noindent The values of the $\gamma_k$ coefficients are obtained by substituting relations (\ref{RelCoeffDRP1}) into (\ref{RelDer}):

\begin{equation}
\label{RelCoeffDRP} \displaystyle \sum _{k=-m}^m \gamma_{k}  =0
\end{equation}

\noindent $m$ being a strictly positive integer, a ${2m+1}$-points
\textit{DRP} scheme (\cite{Tam}) is thus given by:

\footnotesize

\begin{equation}
\label{DRP}
\begin{aligned}
&-u_{i}^{n+1}+u_{i}^{n } +\frac{ \tau }{h}\, \displaystyle \sum
_{k=-m}^m \gamma_{k}\,u_{i+k}^{n }  =0
   \end{aligned}
\end{equation}

\normalsize

 \noindent where the $\gamma_{k}$, $k\in \{-m,m\}$ are the coefficients of the considered
 scheme, and satisfy the relations (\ref{RelCoeffDRP1}).

\section{General study of \textit{DRP} schemes}
\label{CaustiquesDRP}

 \noindent The dispersion
relation related to a general \textit{DRP}-scheme (\ref{DRP}) is
given by:

\footnotesize

\begin{equation}
\begin{aligned}
&\frac{ \tau }{h}\, \displaystyle \sum _{k=-m}^m \gamma_{k}\,e^{i
\,\left(k\, \varphi +\xi _{\omega }\,\tau \right)-B \,\tau +e^{i
\,\xi _{\omega }\,\tau-B \,\tau}}-1=0
\end{aligned}
\end{equation}

\normalsize

 \noindent from which it comes that

\footnotesize
\begin{equation}\begin{aligned}
& i\,{ \xi _{\omega }\,\tau }=B \,\tau -\ln \left(1+\frac{ \tau
\displaystyle \sum _{k=-m}^m \gamma_{k}\,e^{ \,i\,k\,\varphi} }{h}
\right)
   \end{aligned}
\end{equation}

\normalsize

\noindent The group velocity can be expressed as

\footnotesize

\begin{equation}
V_g = \frac{i \,\displaystyle \sum _{k=-m}^m i\,k
\,\gamma_{k}\,e^{\,i\,k\,\varphi  } }{\omega \left(\frac{\sigma
\displaystyle \sum
   _{k=-m}^m \gamma_{k}\,e^{\,i\,k\,\varphi }  }{c}+1\right)}
\end{equation}

\normalsize

 \noindent from which it comes that

\footnotesize
\begin{equation}
\label{rel}
 \begin{aligned} \frac{{\partial } V_g }{\partial
{\varphi}}&=  \frac{i \,c^2 \tau\, \left(\left(c+\sigma
\,\displaystyle  \sum _{k=-m}^m \gamma_{k}\,e^{i\,k\,\varphi }
\right)
  \displaystyle  \sum _{k=-m}^m \, i\,k ^2 \gamma_{k}\,e^{i\,k\,\varphi }  -\sigma \, \left(\displaystyle \sum _{k=-m}^m
   \gamma_{k}\,e^{ i\,k \varphi  }  i\,k  \right){}^2\right)}{\sigma \,   \left(c+\sigma  \,\displaystyle \sum
   _{k=-m}^m \gamma_{k}\,e^{i\,k\,\varphi   } \right){}^2}
   \end{aligned}
\end{equation}

\normalsize

\noindent Through identification of the real and imaginary part of
(\ref{rel}), we obtain:

\footnotesize
\begin{equation}
\label{rel1} \sigma  \,  \displaystyle  \sum _{k,l=-m}^m \,
\gamma_{k}\,\gamma_{i+l}\,\left \lbrace k^2 \sin \left [ \,(k+l)
\varphi \,\right ] -k\, l\, \cos \left [\,(k+l) \varphi\,
   \right ]\right \rbrace =-c \,  \displaystyle  \sum _{k=-m}^m k^2 \,\gamma_{k}\,\sin (k \,\varphi
   )
\end{equation}

\noindent and

\footnotesize
\begin{equation}
\label{rel2}
 \sigma  \, \displaystyle  \sum _{k,l=-m}^m \gamma_{k}\,\gamma_{l}\,   \left \lbrace -\,\cos \left [\,(k+l)\,
  \varphi\right
  ]
    -k\,l\,\sin \left [\,(k+l)\, \varphi \,\right ] \right \rbrace=c \,  \displaystyle  \sum _{k=-m}^m k^2 \,\gamma_{k}\,\cos (k \,\varphi
   )
\end{equation}

\normalsize

\noindent Due to (\ref{RelCoeffDRP1}), (\ref{rel1}) and
(\ref{rel2}) respectively become:

\footnotesize
\begin{equation}
\label{rel1} \sigma  \,  \displaystyle  \sum _{k,l=-m}^m \,
\gamma_{k}\,\gamma_{l}\,\left \lbrace k^2 \sin \left [ \,(k+l)
\varphi \,\right ] -k\, l\, \cos \left [\,(k+l) \varphi\,
   \right ]\right \rbrace =-2\, c \,  \displaystyle  \sum _{k=1}^m k^2 \,\gamma_{k}\,\sin (k \,\varphi
   )
\end{equation}

\noindent and

\footnotesize
\begin{equation}
\label{rel2}
 \sigma  \, \displaystyle  \sum _{k,l=-m}^m \gamma_{k}\,\gamma_{l}\,   \left \lbrace -\,\cos \left [\,(k+l)\,
  \varphi\right
  ]
    -k\,l\,\sin \left [\,(k+l)\, \varphi \,\right ] \right
    \rbrace=0
\end{equation}

\normalsize

\noindent Denote by $T_j$, $j\in  \N^*$ the Chebyshev polynomial
of the first kind, and by $U_j$, $j\in  \N^*$ the Chebyshev
polynomial of the second kind:

\begin{equation}
\label{Tche1} \cos(j\,x)=T_j\left (\cos(x)\right )=\frac
{n}{2}\displaystyle \sum_{k=0}^{\big [ \frac {n}{2}\big ]}
(-1)^k\,\frac {(n-k-1)!}{ k!\,(n-2\,k)!}\,(2\,\cos(x))^{n-2k}
\end{equation}

\begin{equation}
\label{Tche2} \sin(j\,x)=\sin(x)\,U_j\left (\cos(x)\right )
\end{equation}

\noindent where:
\begin{equation}
\label{Tche2} U_j\left (\cos(x)\right )= \displaystyle
\sum_{k=0}^{\big [ \frac {n}{2}\big ]} (-1)^k\,\frac {(n-k )!}{
k!\,(n-2\,k)!}\,(2\,\cos(x))^{n-2k}
\end{equation}

\noindent $\big [ \frac {n}{2}\big ]$ denotes the integer part of
$ \frac {n}{2} $.

\noindent Equations (\ref{rel1}), (\ref{rel2}) can thus be written
as:

\footnotesize
\begin{equation}
\label{rel1bis} \sigma  \,  \displaystyle  \sum _{k,l=-m}^m \,
\gamma_{k}\,\gamma_{l}\,\left \lbrace k^2
\sin(\varphi)\,U_{k+l}\left (\cos(\varphi)\right ) -k\, l\,
T_{k+l}\left (\cos(\varphi)\right )\right \rbrace =-c \,
\displaystyle  \sum _{k=-m}^m k^2
\,\gamma_{k}\,\sin(\varphi)\,U_{k}\left (\cos(\varphi)\right )
\end{equation}

\normalsize

\noindent and

\footnotesize
\begin{equation}
\label{rel2bis}
 \sigma  \, \displaystyle  \sum _{k,l=-m}^m \gamma_{k}\,\gamma_{l}\,   \left \lbrace  T_{k+l}\left (\cos(\varphi)\right )
    +k\,l\,\sin(\varphi)\,U_{k+l}\left (\cos(\varphi)\right )\right
    \rbrace=0
\end{equation}

\normalsize

\noindent Using the relation:

\begin{equation}
\sin(\varphi)=\sqrt{1-\cos^2(\varphi)}
\end{equation}

\noindent equations (\ref{rel1bis}), (\ref{rel2bis}) can be
written as:

\begin{equation}
\label{rel1ter} f_1\left (\cos(\varphi)\right )=0
\end{equation}

\noindent and

\begin{equation}
\label{rel2ter}
 f_2\left (\cos(\varphi)\right )=0
\end{equation}

\noindent where, for all $\theta \in \R$:

\footnotesize
\begin{equation}
f_1(\theta) =\sigma  \,  \displaystyle  \sum _{k,l=-m}^m \,
\gamma_{k}\,\gamma_{ l}\,\left \lbrace k^2
\sqrt{1-\theta^2}\,\,U_{k+l}\left (\theta\right ) -k\, l\,
T_{k+l}\left (\theta\right )\right \rbrace \,+\,c \, \displaystyle
\sum _{k=-m}^m k^2 \,\gamma_{k}\,\sqrt{1-\theta^2}\,U_{k}\left
(\theta\right )
\end{equation}

 \normalsize

\noindent i.e.:
\footnotesize
\begin{equation}
f_1(\theta) =\sigma  \,  \displaystyle  \sum _{k,l=-m}^m \,
\gamma_{k}\,\gamma_{ l}\,\left \lbrace k^2
\sqrt{1-\theta^2}\,\,U_{k+l}\left (\theta\right ) -k\, l\,
T_{k+l}\left (\theta\right )\right \rbrace \,+2\,c \, \displaystyle
\sum _{k=1}^m k^2 \,\gamma_{k}\,\sqrt{1-\theta^2}\,U_{k}\left
(\theta\right )
\end{equation}

 \normalsize

\noindent and

\footnotesize
\begin{equation}
f_2(\theta) =  \displaystyle  \sum _{k,l=-m}^m
\gamma_{k}\,\gamma_{ l}\,   \left \lbrace  T_{k+l}\left
(\theta\right )
    +k\,l\,\sqrt{1-\theta^2}\,U_{k+l}\left (\theta \right )\right \rbrace
\end{equation}

 \normalsize

 \noindent Due to:

\begin{equation}
   T_{j}\left
(1\right )=1 \,\,\,    \forall\, j\in\N^*
\end{equation}
 \noindent it is worth noting that:

\begin{equation}
f_1(1) =-\sigma  \,  \displaystyle  \sum _{k,l=-m}^m \,
\gamma_{k}\,\gamma_{ l}\, k\, l
\end{equation}

\noindent and

\begin{equation}
f_2(1) =  \displaystyle  \sum _{k,l=-m}^m
\gamma_{k}\,\gamma_{ l}
  \end{equation}

 \normalsize

 \noindent The knowledge of the scheme coefficients $\gamma_{k}$, $k\in \{-m,m\}$, enables one to study their variations, and to determine wether the
 equations (\ref{rel1ter}), (\ref{rel2ter}) admit a solution. One
 can thus know wether $\frac{{\partial } V_g }{\partial
{\varphi}}=0$ admits real roots, i. e. wether the schema has
spurious caustics.

\section{Numerical application: the 3-points DRP scheme}
\label{Example} 

\noindent The 3-points DRP scheme is given by:

\begin{equation}
\gamma_1= 0.63662
\end{equation}

 \noindent We thus have:

\begin{equation}
f_1(1) =-2\,\sigma  \, \left \lbrace \gamma_1^2-\gamma_1^2\right \rbrace=0
\end{equation}

\noindent and

\begin{equation}
f_2(1) =  2\, \left \lbrace \gamma_1^2-\gamma_1^2\right \rbrace=0
  \end{equation}

\normalsize

 \noindent For the 3-points DRP scheme, the dispersion
relation is:

\footnotesize

\begin{equation}
\begin{aligned}
&  \,e^{\,i\, \varphi } \left(-e^{-\eta_{\omega} \,\tau}+e^{i
\,\tau\, \xi _{\omega }}\right)+e^{i
   \tau\,\xi _{\omega }} \left(-0.63662 +0.63662 \, e^{2\,i \,\varphi }\right)\, \sigma=0
\end{aligned}
\end{equation}

\normalsize

 \noindent which leads to:

\begin{equation}
\begin{aligned}
&   e^{i
\,\tau\, \xi _{\omega }}    =\frac{e^{\,i\, \varphi } \,e^{-\eta_{\omega} \,\tau}}{ \,e^{\,i\, \varphi }  + 0.63662 \,\sigma\, (e^{2\,i \,\varphi }-1)}
\end{aligned}
\end{equation}

 \noindent It yields:

\begin{equation}
\begin{aligned}
&      \xi _{\omega }     =\frac{1}{\tau}\,\text{Arctan} \left [\frac {
 - \left (1+     0.63662
\,\sigma \, \right ) \,\sin(\varphi ) }
 {1+\left ( 0.63662 \,\sigma -1\right )\,\cos (\varphi )
  } \right ]
 \end{aligned}
\end{equation}

\footnotesize

\normalsize

\noindent The derivative $\frac{{\partial } V_g }{\partial
{\varphi}}$ of the group velocity $V_g$ vanishes for $\varphi=0$, $\varphi=\frac{\pm\pi}{2}$, and $\varphi=0.950935$.

\noindent The 3-points DRP scheme admits thus spurious caustics.

\noindent We now illustrate the caustic phenomenon considering the
two following sinusoidal wave packets:

\begin{equation}
\label{paquet1} u_1=e^{\,-\alpha\,
(x-x_0^1-c\,t)^2}\,\text{Cos}\,[\,k_1\,(x-x_0^1-c\,t)\,]
\end{equation}

\begin{equation}
\label{paquet2} u_2=e^{\,-\alpha\,
(x-x_0^2-c\,t)^2}\,\text{Cos}\,[\,k_2\,(x-x_0^2-c\,t)\,]
\end{equation}

\noindent where $\alpha>0$. The two wave packets are initially
centered  at $x_0^1$ and $x_0^2$, respectively. The group velocity
of the two wave packets are $V_1 = V_g (k_1)$ and  $V_2 = V_g
(k_2)$, respectively, where the function $V_g (x)$ is associated
to the numerical scheme used to solve Eq. (\ref{transp}).

\noindent If the solution obeys the linear advection law given by
Eq. (\ref{transp}), the initial field is passively advected at
speed $c$, while, if the advection speed is scale-dependent (as in
numerical solutions), the two packets will travel at different
speeds, leading to the rise of discrepancies with the
constant-speed solution. Another dispersive error is the
shape-deformation phenomenon: due to numerical errors, the exact
shape of the wave packets will not be exactly preserved. This
secondary effect will not be considered below, since it is not
related to the existence of spurious caustics. It is emphasized
here that the occurance of spurious caustics originates in the
differential error in the group velocity, not in the fact that
shapes of the envelope of the wave packets are not preserved. The
issue of deriving shape-preserving schemes for passive scalar
advection has been adressed by several authors (e.g.
\cite{leonard1,leonard2}).

\noindent The spurious caustics will appear if the two wave
packets happen to get superimposed. During the cross-over, the
$L_\infty$ norm of the error (defined as the difference between
the constant-speed solution and the dispersive one) will exhibit a
maximum. The  characteristic life time of the caustic, $t^*$,
depends directly on the difference between the advection speeds of
the two wave packets and the wave packet widths. Denoting $l_1$
and $l_2$ the characteristic length of the two wave packets, the
time during which they will be (at least partially) superimposed
can be estimated as

\begin{equation}
t^* = \frac{l_1 + l_2}{\vert V_1 - V_2 \vert}
\end{equation}

\noindent It is seen that, since caustics are defined as solutions
for which $\partial V_g / \partial k =0$, $t^*$ will be large if $
\vert k_1 - k_2 \vert \ll 1$. Noting $ k_1 = k_c + \delta k $ and
$ k_2 = k_c - \delta k$, one obtains

\begin{equation}
t^* \simeq \frac{l_1 + l_2}{ 2  ( \delta k) ^2 \left| \frac{\partial
^2 V_g}{\partial k} (k_c) \right| }
\end{equation}

\noindent leading to $t^* \propto  ( \delta k ) ^{-2} $.

\noindent Neglecting shape-deformation effects and assuming that
the numerical scheme is non-dissipative, the numerical error $E$
is given by:

\begin{equation}
\begin{aligned}
E = \vert & \, e^{- \alpha  (x-x_0^1-c\,t )^2} \text{Cos} [
k_1\,(x-x_0^1-ct)\,] -  e^{- \alpha
(x-x_0^1-t \,{V_1})^2} \text{Cos} [ k_1 (x-x_0^1-t {V_1})] \\
&+ e^{- \alpha (x-x_0^2-c\,t)^2} \text{Cos} [ k_2 (x-x_0^2-c t) ] -
e^{ - \alpha (x-x_0^2-t \,{V_2})^2} \,
\text{Cos}\,[\,k_2\,(x-x_0^2-t \,{V_2})\,] \vert
\end{aligned}
\end{equation}

\noindent A simple analysis show that

\begin{equation}
\lim _{t \rightarrow + \infty} L_\infty ( E(t) ) = L_\infty ( u_1
(t=0)) , \quad \max _t  L_\infty ( E(t) ) = 2 L_\infty ( u_1 (t=0))
\end{equation}

\begin{figure}[htbp]
\center{\includegraphics[height=4cm]{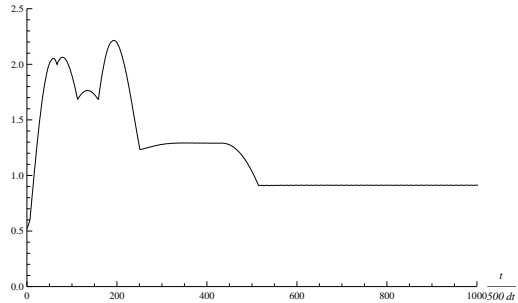} }
\caption{Time history ($L_\infty$ norm) of the numerical error for the two-wave
packet problem (shape deformation and dissipative errors are
neglected to emphasize the linear focusing phenomenon). Numerical parameters are $\alpha
=0.0005$, $h = 0.01$, $V_1 = -2.68381$, $V_2 =-2.51381 $, corresponding to the
properties of the 3-points \textit{DRP} scheme, for $\sigma =0.9$.}
\label{fig-error1}
\end{figure}

\noindent The time history of the $L_\infty$ norm of
$E$ for the 3-points \textit{DRP} scheme scheme,
is displayed in Fig. \ref{fig-error1}, showing the occurance of the
caustic and the sudden growth of the $L_\infty$ error norm.

\begin{figure}[h!]
\center{\includegraphics[height=5cm]{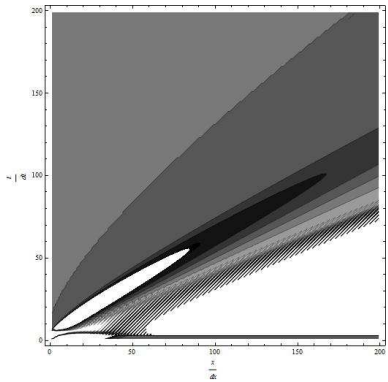}}\\
\caption{\small{Isovalues of the residual kinetic energy for the
3-points \textit{DRP} scheme, for $cfl=0.9$}.} \label{CaustDRP3}
\end{figure}

\noindent Figure \ref{CaustDRP3} displays the isovalues of the
residual kinetic energy for 3-points \textit{DRP} scheme, for $cfl=0.9$. Minima are in black, maxima in white. In
each case, the caustic corresponds to the white domain, where the
residual kinetic energy is maximal.

\section{Concluding remarks}

In the above, we have set a general method that enables
one to determine wether a \textit{DRP} scheme admits or not spurious caustics.\\
The existence of spurious numerical caustics in linear advection \textit{DRP}
schemes has been proved. This linear dispersive phenomenon gives
rise to a sudden growth of the $L_\infty$ norm of the error, which
corresponds to a local focusing of the numerical error in both space
and time. In the present analysis, spurious caustics have been shown
to occur in polychromatic solutions.
The energy of the caustic phenomenon depends on the number of
spectral modes that will get superimposed at the same time. As a
consequence, the spurious error pile-up will be more pronounced  in
simulations with very small wave-number increments.
It has been shown that a popular existing scheme, as the 3-points \textit{DRP}-scheme, allows the
existence of spurious caustics.

\end{document}